\documentclass[a4paper]{article}

\usepackage{amsmath}
\usepackage{amsfonts}
\usepackage{amssymb}         
\usepackage{amsopn}
\usepackage{theorem}          
\usepackage{latexsym}
\usepackage{graphicx}        
\usepackage{mathrsfs}           

\theoremstyle{changebreak}                

\newenvironment{proof}
 {{\sl Proof.}\hspace*{1 ex}}%
 {{\nopagebreak\hspace*{\fill}$\Box$\par\vspace{12pt}}}
\begin{document}

\thispagestyle{empty}
\begin{center} 

{\LARGE A Theorem on Prime Numbers}
\par \bigskip
{\sc Leo Liberti} 
\par  
{\it Centre for Process Systems Engineering
\\
Imperial College of Science, Technology and Medicine}
\par
({\tt l.liberti@ic.ac.uk})
\par \medskip 2 August 2002
\end{center}
\par \bigskip

\begin{abstract} 
The theorem presented in this paper allows the creation of large prime
numbers (of order up to $o(n^2)$ given a table of all primes up to $n$. 
\end{abstract}


\ \hfill \hrule \hfill \ 

Notation: in what follows, products taken over empty index sets are to 
be considered equal to 1.

\noindent
{\bf Theorem} \\
{\sl Let $p(i)$ be the $i$-th prime number and let $I_1,I_2$ be a partition 
of $\{1,\ldots , n\}$ such that}
\begin{eqnarray}
  q_1 = \prod_{i\in I_1} p(i) - \prod_{i\in I_2} p(i) &\le & (p(n))^2, \\
  q_2 = \prod_{i\in I_1} p(i) + \prod_{i\in I_2} p(i) &\le & (p(n))^2.
\end{eqnarray}
{\sl Then $q_1,q_2$ are prime numbers.}

\begin{proof}
Suppose there is a non-unit prime $b\in\mathbb{Z}$ such that $b\le
\sqrt{q_1}$ and $b|q_1$. Then because $\sqrt{q_1}\le p(n)$ we have
$b\le p(n)$; thus there is a $j\le n$ such that $b=p(j)$. Assume
without loss of generality $j\in I_1$ (a symmetric argument holds if
we assume $j\in I_2$). Then $b|q_1$ and $b|\prod_{i\in I_1}p(i)$ imply
$b|\prod_{i\in I_2}p(i)$, i.e. $j\in I_1\cap I_2$, which is empty, so
such a $b$ cannot exist. Hence $q_1$ is prime. Similarly for $q_2$.
\end{proof}

This theorem allows us, given a table of prime numbers up to an
integer $n$, to create prime numbers of order at most $o(n^2)$. 

\noindent {\bf A note about the theorem, by C. Helfgott}

To my posting this theorem on {\tt arXiv.org}, Dr. C. Helfgott of
Berkeley, USA, gave this reply:

\begin{quote}
Your results, while not wrong, are virtually useless. I do not mean to
offend, but rather to point out certain facts you may have
inadvertantly overlooked:
\begin{enumerate} 
\item You define a quantity
$$
 q_2 = \Pi_{i \in I_1} p(i) + \Pi_{i \in I_2} p(i)
$$
and wish to restrict it to the range $\{0, \ldots, (p(n))^2\}$.
However, for any value of $n >5$, and any partition $I_1, I_2$ of
$\{1,\ldots n\}$, the quantity $q_2$ will be greater than your
bound. Thereby making your theorem a null statement.
\item
For your quantity $q_1$, yes, it may be possible to construct a partition
with the properties described. I suspect, however, that you would find it
computationally intractable to do so for any significant value of $n$.
Furthermore, there are much better ways of constructing primes of size
$p(n)^2$ than to start by listing all primes up to $p(n)$ (which is in and of
itself a nigh-impossible task for any significant value of $n$).
\end{enumerate}
In conclusion, while the two results presented in your paper are true,
one is a null statement, and the other is so impractical to implement as
to be completely useless.
\end{quote}

I can only thank Dr. Helfgott for highlighting these facts. I would
just like to point out that since the inception of the polynomial-time
primality testing algorithm by Agrawal, Kayal and Saxena (see {\tt
http://www.cse.iitk.ac.in/news/primality.html}), constructing a list
of all primes up to $p(n)$ is less ``nigh-impossible'' than it
seemed. It is true, however, that the few numerical experiments I ran
on this theorem were disappointing.


\bibliographystyle{alpha}
\bibliography{phd}

\end{document}